\documentclass[11pt,notitlepage,twoside]{article}
\usepackage{latexsym}
\usepackage{amsmath}
\usepackage{amsfonts}
\usepackage{amssymb}
\renewcommand{\a }{\alpha }
\renewcommand{\b }{\beta }
\renewcommand{\d}{\delta }

\newcommand{\D }{\Delta }

\newcommand{\e }{\varepsilon }
\newcommand{\g }{\gamma}

\renewcommand{\l }{\lambda }

\newcommand{\n }{\nabla }
\newcommand{\var }{\varphi }

\renewcommand{\th }{\theta }
\renewcommand{\o }{\omega }
\renewcommand{\O }{\Omega }

\newcommand{\ov}{\overline}

\newcommand{\be}{\begin{equation}}
\newcommand{\ee}{\end{equation}}
\newenvironment{pf}{\noindent{\bf Proof.}\enspace}{
\hfill$\Box$\medskip}
\newenvironment{pfn}[1]{\noindent{\bf Proof of {#1}\enspace}}{
\hfill$\Box$\medskip}
\newcommand{\R}{\mathbb{R}}
 
\newtheorem{thm}{Theorem}[section]
\newtheorem{pro}[thm]{Proposition}
\newtheorem{lem}[thm]{Lemma}

\numberwithin{equation}{section}

\author{{\large \textbf{Khalil EL MEHDI} }\\
{\it\small Facult\'e des Sciences et Techniques}\\
 {\it\small Universit\'e de Nouakchott, BP 5026}\\
 {\it\small Nouakchott, Mauritania}\\
{\it\small E-mail: khalil@univ-nkc.mr}\\
{\it\small and }\\
{\it\small The Abdus Salam ICTP, Trieste, Italy }
}

\title { {\Large \textbf{Single Blow-up Solutions for a Slightly Subcritical\\ Biharmonic Equation  }} }

\begin{document}
\date{ }
\maketitle
{\footnotesize
\noindent{\bf Abstract.}
In this paper, we consider a biharmonic equation under the Navier boundary condition and with a nearly critical exponent  $(P_{\e})$:  $\D^2u=u^{9-\e}$, $u>0$ in $\O$ and $u=\D u=0$ on $\partial\O$, where $\O$ is a smooth bounded domain in $\R^5$, $\e >0$. We study the asymptotic behavior of solutions of $(P_{\e})$ which are minimizing for the Sobolev quotient as $\e$ goes to zero. We show that such solutions concentrate around a point $x_0 \in \O$ as $\e\to 0$, moreover $x_0$ is a critical point of the Robin's function. Conversely, we show that for any nondegenerate critical point $x_0$ of the Robin's function, there exist solutions of $(P_{\e})$ concentrating around $x_0$ as $\e\to 0$. 
}

\section{Introduction and   Results }
\mbox{}
Let us consider the following  biharmonic equation under the Navier boundary condition
$$
(Q_{\e} ) \quad \left\{
\begin{array}{cc}
 \D ^2 u = u^{p-\e},\,\, u>0 &\mbox{ in }\, \O \\
    \D u= u =0   \quad\quad     & \mbox{ on }\,  \partial  \Omega ,
\end{array}
\right.
$$
where $\O$ is a smooth bounded domain in $\R^n$, $n\geq 5$, $\e$ is a small positive parameter, and $p+1= 2n/(n-4)$ is the critical Sobolev exponent of the embedding  $H^2(\O )\cap H^1_0(\O ) \hookrightarrow L^{2n/(n-4)}(\O )$.

It is known that $(Q_\e)$ is related to the limiting problem $(Q_0)$ (when $\e=0$) which exhibits a lack of compactness and gives rise to solutions of $(Q_\e)$ which blow up as $\e\to 0$. The interest of the limiting problem $(Q_0)$ grew  from its resemblance to some geometric equations involving  Paneitz operator and which has widely been studied in these last years (for details one can see \cite{BE}, \cite{BE3}, \cite{C}, \cite{DHL}, \cite{DMA1}, \cite{DMA2}, \cite{F} and  references therein).

Several authors have studied the existence and behavior of blowing up solutions for the corresponding second order elliptic problem (see, for example, \cite{AP}, \cite{AA}, \cite{BLR}, \cite{BP}, \cite{H}, \cite{L}, \cite{R1}, \cite{R2}, \cite{R3} and  references therein). In sharp contrast to this, very little is known for fourth order elliptic equations. In this paper we are mainly interested in the asymptotic behavior and the existence of solutions of $(Q_\e)$ which blow up around one point, and the location of this blow up point as $\e\to 0$.

The existence of solutions of $(Q_\e)$ for all $\e\in (0,p-1)$ is well known for any domain $\O$ (see, for example \cite{EFT}). For 
  $\e =0$, the situation is more complex, Van Der Vorst showed in \cite{V1} that if $\O$ is starshaped $(Q_{0})$ has no solution whereas Ebobisse and Ould Ahmedou proved in \cite{EO} that  $(Q_{0})$ has a solution provided that some homology group of $\O$ is nontrivial. This topological condition is sufficient, but not necessary, as examples of contractible domains $\O$ on which a solution exists show \cite{GGS}.

In view of this qualitative change in the situation when $\e=0$, it is interesting to study the asymptotic behavior of the subcritical solution $u_\e$ of  $(Q_{\e})$ as $\e\to 0$. Chou-Geng \cite{CG}, and Geng \cite{G} made a first study, when $\O$ is strictly convex. The convexity assumption was needed in their proof in order to apply the method of moving planes (MMP for short) in proving a priori estimate near the boundary. Notice that in the Laplacian case (see \cite{H}), the MMP has been used to show that blow up points are away from the boundary of the domain. The process is standard if domains are convex. For nonconvex regions, the MMP still works in the Laplacian case through the applications of Kelvin transformations \cite{H}. For $(Q_{\e})$, the MMP also works for convex domains \cite{CG}. However, for nonconvex domains, a Kelvin transformation does not work for  $(Q_{\e})$ because the Navier boundary condition is not invariant under the Kelvin transformation of biharmonic operator. In \cite{BE2}, Ben Ayed and El Mehdi removed the convexity assumption of Chou and Geng for higher dimensions, that is $n\geq 6$. The aim of this paper is to prove that the results of \cite{BE2} are true in dimension $5$. In order to state precisely our results, we need to introduce some notations.

We consider the following problem
$$
(P_{\e} ) \quad \left\{
\begin{array}{cc}
 \D ^2 u = u^{9-\e},\,\, u>0 &\mbox{ in }\, \O \\
    \D u= u =0   \quad\quad     & \mbox{ on }\,  \partial  \Omega ,
\end{array}
\right.
$$
where $\O$ is a smooth bounded domain in $\R^5$ and $\e$ is a small positive parameter.

Let us define on $\O$ the following Robin's function
$$
\var (x)=H(x,x),\quad \mbox{with}\quad H(x,y)=|x-y|^{-1} -G(x,y),\,\,\mbox{for }\,\, (x,y)\in\O\times\O,
$$
where $G$ is the Green's function of $\D^2$, that is,
$$
\forall x\in\O \quad \left\{
\begin{array}{cc}
 \D ^2 G(x,.) = c\d_x &\mbox{ in }\, \O \\
    \D G(x,.)= G(x,.) =0      & \mbox{ on }\,  \partial  \Omega ,
\end{array}
\right.
$$
where $\d_x$ denotes the Dirac mass at $x$ and $c=3\o_5$, with $\o_5$ is the area of the unit sphere of $\R^5$.\\
For $\l>0$ and $a\in\R^5$, let
\begin{eqnarray}\label{e:11}
 \d _ {a,\l }(x) =  \frac {c_0\lambda
 ^{1/2}}{(1+\lambda^2|x-a|^2)^{1/2}},\,\,
 c_0=(105)^{1/8}.
\end{eqnarray}
It is well known (see \cite{Lin}) that $\d_{a,\l}$ are the only
solutions of
\begin{eqnarray*}
 \D^2 u =  u^{9},\quad  u>0 \mbox{  in  } \R^5
\end{eqnarray*}
\noindent
and are also the only minimizers of the Sobolev inequality on the
whole space, that is
\begin{eqnarray}\label{e:12}
 S =\inf\{|\D u|^{2}_{L^2(\R^5)}|u|^{-2}_{L^{10}(\R^5)}
,\, s.t.\, \D u\in L^2 ,u\in L^{10} ,u\neq 0 \}.
\end{eqnarray}
\noindent
 We denote by  $P\d _{a,\l}$ the projection of $\d
_{a,\l}$ on $\mathcal{H}(\O):=H^2(\O )\cap H^1_0(\O)$, defined by
$$ \D^2 P\d_{a,\l}=\D^2\d_{a,\l} \mbox{ in } \O  \mbox{ and }
\D P\d_{a,\l}=P\d_{a,\l}=0\mbox{ on } \partial \O.
$$
Let
\begin{align}
\th_{a,\l} &= \d_{a,\l} -P\d_{a,\l},\label{th}\\
||u||&=\left(\int_\O |\D u|^2\right)^{1/2},\quad \langle u,v\rangle =\int_\O \D u\D v,\quad u,v \in H^2(\O )\cap H^1_0(\O)  \label{e:13}\\
||u||_q&=|u|_{L^q(\O )} \label{e:15}.
\end{align}

Thus we have the following result:
\begin{thm}\label{t:11} 
Let $(u_\e)$ be a solution of $(P_{\e})$, and assume that
$$
{|| u_\e||^2}{||u_\e||_{10-\e}^{-2}}\to S \mbox{ as } \e \to
0,\leqno{(H)}
$$
where $S$ is the best Sobolev constant in $\R^5$ defined by \eqref{e:12}. Then (up to a
subsequence) there exist $a_\e\in \O$, $\l_\e > 0$, $\a_\e > 0$
and $v_\e$ such that $u_\e$ can be written as
$$
u_\e = \a_\e P\d _{a_\e , \l _\e}+v_\e
$$
with $\a_\e\to 1$, $||v_\e|| \to 0$, $a_\e \in
\O$ and $\l _\e d(a_\e, \partial\O ) \to +\infty$ as $\e \to 0$.\\
In addition, $a_\e$ converges to a critical point $x_0\in \O$ of
$\var$ and we have
$$
\lim_{\e\to 0}\e||u_\e||_{L^\infty (\O )}^2 = (c_1c_0^2/c_2)\var
(x_0),
$$
 where  $c_1=c_0^{10} \int_{\R^5} \frac{d x
}{(1+|x|^2)^{9/2}}$, $c_2=c_0^{10}\int_{\R^5}
\frac{\log(1+|x|^2) (1-|x|^2)}{(1+|x|^2)^{6}} dx $ and $c_0=(105)^{1/8}$.
\end{thm}

Our next result provides a kind of converse to Theorem \ref{t:11}.
\begin{thm}\label{t:12} Assume that  $x_0\in\O$ is
a nondegenerate critical point of $\var$. Then
there exists an $\e_0 >0$ such that for each $\e\in (0,\e_0]$,
$(P_{\e})$ has a solution of the form
$$
u_\e = \a_\e P\d _{a_\e , \l _\e}+v_\e
$$
with $\a_\e\to 1$, $||v_\e|| \to 0$, $a_\e \to
x_0$ and $\l _\e d(a_\e, \partial\O ) \to +\infty$ as $\e \to 0$.
\end{thm}

Our strategy to prove the above results is the same as in higher dimensions. However, as usual in elliptic equations involving critical Sobolev exponent, we need more refined estimates of the asymptotic profiles of solutions when $\e \to 0$ to treat the lower dimensional case. Such refined estimates, which are of self interest, are highly nontrivial and use in a crucial way  careful expansions of the Euler-Lagrange functional associated to $(P_\e)$, and its gradient near a small neighborhood of highly concentrated functions. To perform such expansions we  make use of the techniques developed by Bahri \cite{B} and Rey \cite{R1}, \cite{R4} in the framework of the {\it Theory of critical points at infinity}.

The outline of the paper is the following: in Section 2 we  perform some crucial estimates needed in our proofs and Section 3 is devoted to the proof of our  results.

\section{Some Crucial Estimates}
\mbox{}
In this section, we prove some crucial estimates which will play an important role in proving our results. We first recall some results.
\begin{pro}\label{p:21}\cite{BH}
Let $a\in \O$ and $\l>0$ such that $\l d(a,\partial \O)$ is large
enough. For $\th_{(a,\l)}=\d_{(a,\l)}-P\d_{(a,\l)}$, we have the
following estimates
$$
(a)\quad 0\leq \th_{(a,\l)}\leq\d_{(a,\l)}, \qquad (b) \quad
\th_{(a,\l)}=c_0 \l^{\frac{-1}{2}}H(a,.) +f_{(a,\l)},$$ where
  $f_{(a,\l)}$ satisfies
$$
f_{(a,\l)}=O\left(\frac{1}{\l^{5/2}d^{3}}\right),\quad
\l\frac{\partial f_{(a,\l)}}{\partial\l}=O\biggl(\frac{1}{
\l^{5/2}d^{3}}\biggr),\quad \frac{1}{\l}\frac{\partial
f_{(a,\l)}}{\partial a}=O\biggl(\frac{1}{
\l^{7/2}d^{4}}\biggr),
$$
 where $d$ is the distance
$d(a,\partial \O)$,
$$
\mid\th_{(a,\l)}\mid_{L^{10}}=O\bigl((\l
d)^{-1/2}\bigr), \quad
\mid\mid\th_{(a,\l)}\mid\mid=O\bigl((\l
d)^{-1/2}\bigr), \leqno{(c)}
$$
$$
 \bigg|\l\frac{\partial\th_{(a,\l)}}{\partial\l}\bigg|_{L^{10}}=O\left(
\frac{1}{(\l d)^{1/2}}\right), \quad
 \bigg|\frac{1}{\l}\frac{\partial\th_{(a,\l)}}{\partial a}
\bigg|_{L^{10}}=O\left(\frac{1}{(\l d)^{3/2}}\right).
$$
\end{pro}
\begin{pro}\label{p:24}\cite{BE2}
Let $u_\e$ be a solution of $(P_{\e})$ which satisfies $(H)$.
Then, there exist $a_\e\in \O$, $\a_\e > 0 $, $\l_\e >0$ and
$v_\e$ such that
$$
u_\e = \a_\e P\d _{a_\e , \l _\e}+ v_\e
$$
with $\a_\e \to 1$, $\l_\e d(a_\e,\partial
\O) \to \infty$, $c_0^{-2}||u_\e||^2_\infty/\l_\e \to 1$, $||u_\e||_\infty ^\e \to 1$ and $||v_\e||\to 0$.\\
Furthermore, $v_\e\in E_{(a_\e,\l_\e)}$ which is the set of  $v\in \mathcal{H}(\O)$ such that  
$$
 \langle  v,P\d_{a_\e,\l_\e}\rangle  =\langle v,\partial P\d_{a_\e, \l_\e}/\partial \l_\e\rangle  =0,\, \langle  v_\e,\partial P\d_{a_\e, \l_\e}/\partial a\rangle  =0 .\leqno{(V_0)}
$$
\end{pro}
\begin{lem}\label{l:23}\cite{BE2}
$\l_\e^\e =1+o(1)$ as $\e$ goes to zero implies that
$$
\d_\e ^{-\e}-c_0^{-\e}\l_\e ^{\e(4-n)/2}=O\left( \e
\log(1+\l_\e ^2 |x-a_\e|^2)\right) \quad \mbox{in}\quad \O.
$$
\end{lem}
\begin{pro}\label{p:26}\cite{BE2}
Let $(u_\e)$ be  a solution of $(P_{\e})$ which satisfies $(H)$.
Then $v_\e$ occurring in Proposition \ref{p:24} satisfies
\be\label{v}
|| v_\e ||\leq C \left(\e +  (\l _\e
    d_\e )^{-1}\right),
\ee
where $C$ is a positive constant independent of $\e$.
\end{pro}

Now, we are going to state and prove the crucial estimates needed in the proof of our theorems.
\begin{lem}\label{l:25}
For $\e$ small, we have the following estimates 
\begin{align*}
i) & \quad \int _{\O} \d _\e ^9 \frac{1}{\l _\e} \frac{\partial P\d
_\e}{\partial a}=-\frac{c_1}{2\l_\e^2} \frac{\partial H}{\partial a}(a_\e,a_\e) + O\left(\frac{1}{(\l_\e d_\e)^3}\right),\\
ii) & \quad \int _{\O} P\d _\e ^{9-\e} \frac{1}{\l _\e} \frac{\partial P\d
_\e}{\partial a}=-\frac{c_1}{\l_\e^{2+\e/2}} \frac{\partial H}{\partial a}(a_\e,a_\e) + O\left(\frac{1}{(\l_\e d_\e)^3}+\frac{\e}{(\l_\e d_\e )^2}\right),
\end{align*}
where $c_1$ is the constant defined in Theorem \ref{t:11}.
\end{lem}
\begin{pf}
Notice that 
\be \label{k:1}
\int_{\O\setminus B_\e} \d_\e ^{10} = O\left(\frac{1}{(\l_\e d_\e)^5}\right).
\ee
Thus, we have, for $1\leq k\leq 5$
\begin{align}\label{1}
 \int _{\O} \d _\e ^9 \frac{1}{\l _\e} \frac{\partial P\d
_\e}{\partial a_k}&=\int _{\O} \d _\e ^9 \frac{1}{\l _\e} \frac{\partial  \d
_\e}{\partial a_k} - \int _{\O} \d _\e ^9 \frac{1}{\l _\e} \frac{\partial  \th
_\e}{\partial a_k}\notag\\
& = - \int _{B_\e} \d _\e ^9 \frac{1}{\l _\e} \frac{\partial  \th
_\e}{\partial a_k} +  O\left(\frac{1}{(\l_\e d_\e)^5}\right),
\end{align}
where $B_\e=B(a_\e,d_\e)$. Expanding $\partial  \th
_\e / \partial a_k$ around $a_\e$ and using Proposition \ref{p:21}, we obtain
\be\label{2}
\int _{B_\e} \d _\e ^9 \frac{1}{\l _\e} \frac{\partial  \th
_\e}{\partial a_k}= \frac{c_0}{2\l_\e^{3/2}}\frac{\partial H(a_\e,a_\e) }{\partial a}\int_{B_\e}\d_\e^9 +  O\left(\frac{1}{(\l_\e d_\e)^3}\right).
\ee
Estimating the integral on the right-hand side in \eqref{2} and using \eqref{1}, we easily derive claim $i)$.\\
To prove claim $ii)$, we write
\begin{align}\label{3}
\int _{\O} P\d _\e ^{9-\e}\frac{1}{\l _\e}& \frac{\partial P\d _\e}{\partial
a_k}= \int _{\O}\d _\e ^{9-\e}\frac{1}{\l _\e} \frac{\partial \d _\e}{\partial a_k} - \int _{\O}\d _\e ^{9-\e}\frac{1}{\l _\e} \frac{\partial \th
_\e}{\partial a_k} -(9-\e)\int _{\O}\d _\e ^{8-\e}\th _\e \frac{1}{\l}
\frac{\partial \d _\e}{\partial a_k} \notag \\
& + \frac{(9-\e)(8-\e)}{2}\int _{\O}\d _\e ^{7-\e}\th _\e^2 \frac{1}{\l} \frac{\partial \d _\e}{\partial a_k}+ O\left(\int _{\O} \d _\e ^{8-\e}\th _\e \bigg|\frac{1}{\l _\e}\frac{\partial \th _\e}{\partial a_k}\bigg|+\int \d_\e ^{7-\e}\th_\e ^3\right)
\end{align}
and we have to estimate each term on the right hand-side of
\eqref{3}.\\
Using Proposition \ref{p:21} and Lemma \ref{l:23}, we have 
\be\label{4}
\int _{\O} \d _\e ^{7-\e}\th _\e^3 \leq c||\th_\e||_\infty ^3 \int \d_\e ^7 =  O\left(\frac{1}{(\l_\e d_\e)^3}\right),
\ee
\be\label{5}
\int _{\O} \d _\e ^{8-\e}\th _\e \bigg|\frac{1}{\l _\e}
\frac{\partial \th _\e}{\partial a_k}\bigg|\leq c ||\th_\e||_\infty \bigg|\bigg|\frac{1}{\l_\e}\frac{\partial \th _\e}{\partial a_k}\bigg|\bigg|_\infty \int_\O \d_\e^8 =  O\left(\frac{1}{(\l_\e d_\e)^3}\right).
\ee
We also have
\be\label{6}
 \int _{\O}\d _\e ^{9-\e}\frac{1}{\l _\e} \frac{\partial \d _\e}{\partial
a_k}=\int _{\O\setminus B_\e}\d _\e ^{9-\e}\frac{1}{\l _\e} \frac{\partial \d _\e}{\partial a_k}=  O\left(\frac{1}{(\l_\e d_\e)^5}\right).
\ee
Expanding $\th_\e$ around $a_\e$ and using Proposition \ref{p:21} and Lemma \ref{l:23}, we obtain
\be\label{7}
9\int _{B_\e}\d _\e ^{8-\e}\th _\e \frac{1}{\l_\e}
\frac{\partial \d _\e}{\partial a_k} = \frac{c_1}{2\l_\e^{2+\e/2}} \frac{\partial H(a_\e,a_\e) }{\partial a} + O\left(\frac{1}{(\l_\e d_\e)^3}+\frac{\e}{(\l_\e d_\e)^2} \right),
\ee
\be\label{m:4}
\int _{B_\e}\d _\e ^{7-\e}\th _\e^2 \frac{1}{\l_\e}
\frac{\partial \d _\e}{\partial a_k} = O\left(\frac{1}{(\l_\e d_\e)^3} \right).
\ee
In the same way, we find
\be\label{8}
\int _{\O}\d _\e ^{9-\e} \frac{1}{\l_\e}
\frac{\partial \th _\e}{\partial a_k} = \frac{c_1}{2\l_\e^{2+\e/2}} \frac{\partial H(a_\e,a_\e) }{\partial a} + O\left(\frac{1}{(\l_\e d_\e)^3}+\frac{\e}{(\l_\e d_\e)^2} \right).
\ee
Combining \eqref{3}--\eqref{8}, we obtain claim $ii)$.
\end{pf}

To improve the estimates of the integrals involving $v_\e$, we use
an  idea of Rey \cite{R4}, namely we write
$$
 v_\e = \Pi v_\e + w_\e, 
$$
 where $\Pi v_\e$
denotes the projection of $v_\e$ onto $H^2 \cap H^1_0(B_\e)$, that is
\be\label{9}
 \D^2 \Pi v_\e=\D^2 v_\e\quad\mbox{in}\quad B_\e;\quad
\D\Pi v_\e=\Pi v_\e= 0\quad\mbox{on}\quad \partial B_\e, 
\ee
 where $B_\e= B(a_\e,d_\e)$. We split $\Pi v_\e$ in an even part
$\Pi v_\e^{e}$ and an odd part $\Pi v_\e^{o}$ with respect to
$(x-a_\e)_k$, thus we have 
\be\label{10}
 v_\e =
\Pi v_\e^{e}+\Pi v_\e^{o} + w_\e \quad\mbox{in}\quad
B_\e\quad\mbox{with}\quad \D^2 w_\e=0\,\,\mbox{in}\,\, B_\e.
 \ee
Notice that it is difficult to improve the estimate \eqref{v} of the $v_\e$-part of solutions. However, it is sufficient to improve the integrals involving the odd part of $v_\e$ with respect to $(x-a_\e)_k$, for $1\leq k\leq 5$ and to know the exact contribution of the integrals containing the $w_\e$-part of $v_\e$. Let us start by the terms involving $w_\e$.
\begin{lem}\label{l:26} For $\e$ small, we have that
$$
\int_{B_\e} \d_\e^8\left(\d_\e^{-\e}-\frac{1}{c_0^\e\l_\e^{\e/2}}\right)\frac{1}{\l_\e}\frac{\partial\d_\e}{\partial a_k} w_\e= O\left(\frac{\e ||v_\e||}{(\l_\e d_\e)^{1/2}}\right).
$$
\end{lem}
\begin{pf}
Let $\psi$ be the solution of
$$
\D^2\psi= \d_\e^8\left(\d_\e^{-\e}-\frac{1}{c_0^\e\l_\e^{\e/2}}\right)\frac{1}{\l_\e}\frac{\partial\d_\e}{\partial a_k} \quad\mbox{in}\quad B_\e; \quad \D \psi=\psi=0\quad\mbox{on}\quad \partial B_\e.
$$
Thus we have 
\be\label{11}
I_\e:= \int_{B_\e} \D^2\psi w_\e= \int_{\partial B_\e}\frac{\partial\psi}{\partial\nu} \D w_\e + \int_{\partial B_\e}\frac{\partial\D\psi}{\partial\nu}w_\e.
\ee
 Let $G_\e$ be the Green's function for the biharmonic operator on $B_\e$ with the Navier boundary conditions, that is,
\be\label{G}
 \D^2 G_\e(x,.)= c\d_x\quad\mbox{in}\quad B_\e;\quad \D G_\e(x,.)=G(x,.)=0\quad\mbox{on}\quad\partial B_\e,
\ee
where $c=3w_5$. Therefore $\psi$ is given by
 $$
\psi(y)=\int_{B_\e}G_\e(x,y)\d_i^8\left(\d_\e^{-\e}-\frac{1}{c_0^\e\l_\e^{\e/2}}\right)\frac{1}{\l_\e}\frac{\partial\d_\e}{\partial a_k},
\quad y\in B_\e
$$
and its normal derivative by
 \be\label{12}
\frac{\partial\psi}{\partial\nu}(y)=\int_{B_\e}\frac{\partial
G_\e}{\partial\nu}(x,y) \d_i^8\left(\d_\e^{-\e}-\frac{1}{c_0^\e\l_\e^{\e/2}}\right)\frac{1}{\l_\e}\frac{\partial\d_\e}{\partial a_k}, \quad y\in \partial B_\e.
 \ee
Notice that:
\begin{align}
&\mbox{for}\,\, x\in B_\e\setminus B(y,d_\e/2),\mbox{ we have } \frac{\partial G_\e}{\partial\nu}(x,y)= O\left(\frac{1}{d_i^2}\right);\qquad  \frac{\partial \D G_\e}{\partial\nu}(x,y)= O\left(\frac{1}{d_i^4}\right) \label{D1}\\
&\mbox{for}\,\, x\in B_\e\cap B(y,{d_\e/2}),\mbox{ we have } \bigg|\frac{\partial G_\e}{\partial\nu}(x,y)\bigg|\leq \frac{c}{|x-y|^2}; \quad  \bigg|\frac{\partial \D G_\e}{\partial\nu}(x,y)\bigg|\leq \frac{c}{|x-y|^4}\label{D2}\\
&\mbox{for}\,\, x\in B_\e\cap B(y, d_\e/2),\mbox{ we
have } \d_i^8\left(\d_\e^{-\e}-\frac{1}{c_0^\e\l_\e^{\e/2}}\right)\frac{1}{\l_\e}\frac{\partial\d_\e}{\partial a_k} = O\left(\frac{\e \log\l_\e d_\e}{(\l_\e d_\e)^9}\right),\notag\\
&\mbox{for}\,\, x\in B_\e\setminus B(y,{d_\e/2}),\mbox{ we have } \d_i^8\left(\d_\e^{-\e}-\frac{1}{c_0^\e\l_\e^{\e/2}}\right)\frac{1}{\l_\e}\frac{\partial\d_\e}{\partial a_k} = O\left(\d_\e^9\e\log(1+\l_\e^2|x-a_\e|^2\right).\notag
\end{align}
Therefore
\be\label{13}
\biggl|\frac{\partial\psi}{\partial\nu}(y)\biggr|=O\left(\frac{\e}{\l_\e^{1/2}d_\e^2}\right).
\ee
In the same way, we have
\be\label{14}
\biggl|\frac{\partial\D\psi}{\partial\nu}(y)\biggr|=O\left(\frac{\e}{\l_\e^{1/2}d_\e^4}\right).
\ee
Using \eqref{11}, \eqref{13}, \eqref{14}, we obtain
\be\label{15}
I_\e=O\left(\frac{\e}{\l_\e^{1/2}d_\e^2}\int_{\partial B_\e}|\D w_\e| + \frac{\e}{\l_\e^{1/2}d_\e^4}\int_{\partial B_\e}| w_\e|\right).
\ee
 To estimate the right-hand side of \eqref{15}, we introduce
the following function
$$
\bar{w}(X)=d_\e ^{1/2}w_\e(a_\e + d_\e X), \quad \bar{v}(X)=d_\e ^{1/2}v_\e(a_\e + d_\e X)\quad \mbox{for } X\in B(0,1).
$$
$\bar{w}$ satisfies 
$$
 \D^2 \bar{w}=0
\quad\mbox{in}\quad B:=B(0,1); \quad
\D\bar{w}=\D \bar{v}, \quad \bar{w}=\bar{v} \quad\mbox{on}\quad \partial B.
$$
 We deduce that 
\be \label{m:1}
 \int_{\partial B}|\D\bar{w}| +  \int_{\partial B}|\D\bar{w}| \leq C\left(\int_B |\D \bar{v}|^2\right)^{1/2}=C\left(\int_{B_\e} |\D v_\e|^2\right)^{1/2} .
\ee
 But, we have 
\be\label{16}
\int_{\partial B}|\D\bar{w}| +  \int_{\partial B}|\D\bar{w}|= 
 \left(\frac{1}{ d_\e}\right)^{3/2} \int_{\partial B_\e}|\D w_\e| + \left(\frac{1}{ d_\e}\right)^{7/2} \int_{\partial B_\e}| w_\e|.
\ee
 Using \eqref{15}, \eqref{m:1} and \eqref{16}, the lemma follows. 
\end{pf}
\begin{lem}\label{l:27}
For 
$\e$ small, we have
\begin{align*}
i)&\quad \int_{B_\e}\D\left(\frac{1}{\l_\e}\frac{\partial \Pi \d_\e}{\partial a_k}\right)\D w_\e= O\left(\frac{||v_\e||}{(\l_\e d_\e)^{3/2}}\right),\\
ii)&\quad \int_{B_\e} \d_\e^{8-\e} \Pi v_\e^o w_\e =  O\left(\frac{||v_\e|| ||\Pi v_\e^o||}{(\l_\e d_\e)^{1/2}}\right).
\end{align*}
\end{lem}
\begin{pf}
Using \eqref{10}, we obtain
\be\label{27i)}
\int_{B_\e}\D\left(\frac{1}{\l_\e}\frac{\partial \Pi \d_\e}{\partial a_k}\right)\D w_\e=\int_{\partial B_\e} \frac{\partial \psi_k}{\partial\nu} \D w_\e,\quad\mbox{with}\quad \psi_k=\frac{1}{\l_\e}\frac{\partial \Pi \d_\e}{\partial a_k}.
\ee
Using an integral representation for $\psi_k$ as in \eqref{12}, we obtain for $y\in \partial B_\e$,
$$
\frac{\partial\psi}{\partial\nu}(y)=\int_{B_\e}\frac{\partial
G_\e}{\partial\nu}(x,y) \D^2 \psi_k,
$$
where $G_\e$ is the Green's function defined in \eqref{G}. Clearly, we have
\be\label{m:2}
\Pi\d_\e(x)=\d_\e(x)-\frac{c_0\l_\e^{1/2}}{(1+\l_\e^2d_\e^2)^{1/2}}- \frac{c_\e(a_\e,d_\e)}{10}(|x-a_\e|^2-d_\e^2),
\ee
with $c_\e(a_\e,d_\e)=\D\d_{\e_{\mid \partial B_\e}}$. Thus we deduce that
\be\label{17}
\frac{\partial\psi}{\partial\nu}(y)=9\int_{B_\e}\frac{\partial
G_\e}{\partial\nu}(x,y)\d_\e^8\frac{1}{\l_\e}\frac{\partial\d_\e}{\partial a_k}.
\ee
In $B_\e\setminus B(a_\e,d_\e/2)$, we argue as in \eqref{13} and \eqref{D2}, we obtain
$$
 \int_{B_\e}\frac{\partial
G_\e}{\partial\nu}(x,y)\d_\e^8\frac{1}{\l_\e}\frac{\partial\d_\e}{\partial a_k} = O\left(\frac{1}{\l_\e^{9/2}d_\e^6}\right).
$$
Furthermore, since
$$
\biggl| \n \frac{\partial
G_\e}{\partial\nu}(x,y)\biggr|=O\left(\frac{1}{d_\e^3}\right)\quad\mbox{for}\quad
(x,y)\in B(a_\e, d_\e/2)\times \partial B_\e,
$$
we obtain
$$
\bigg|\int_{ B(a_\e, d_\e/2)}\frac{\partial G_\e}{\partial\nu} (x,y)\d_\e^8
\frac{1}{\l_\e}\frac{ \partial\d_\e}{\partial a_k}\bigg|
 \leq \frac{c}{d_\e ^3}
 \int_{ B(a_\e,{ d_\e}/{2})} \d_\e
^9 |x-a_\e|= O\left(\frac{1}{\l_\e^{{3}/{2}}d_\e^3}\right),
$$
where we have used the evenness of $\d_\e$ and the oddness of its
derivative. Thus
\be\label{18}
\frac{\partial\psi_k}{\partial
\nu}(y)=O\left(\frac{1}{\l_\e^{3/2}d_\e^3}\right).
 \ee
Using \eqref{27i)} and \eqref{18}, we obtain
$$
\int_{B_\e}\D\left(\frac{1}{\l_\e}\frac{\partial \Pi \d_\e}{\partial a_k}\right)\D w_\e\leq \frac{c}{\l_\e^{3/2}d_\e^3}\int_{\partial B_\e} |\D w_\e|.
$$
Arguing as in \eqref{16}, claim $i)$ follows. To prove claim $ii)$, let $\psi$ be such that
$$
\D^2\psi=\d_\e^{8-\e}\Pi v_\e^o\quad\mbox{in}\quad B_\e;\quad\D\psi=\psi=0\quad\mbox{on}\quad\partial B_\e.
$$
We have
\be\label{19}
\int_{B_\e} \d_\e^{8-\e}\Pi v_\e^o w_\e = \int_{\partial B_\e}\frac{\partial \D\psi}{\partial\nu} w_\e + \int_{\partial B_\e}\frac{\partial\psi}{\partial\nu}\D w_\e.
\ee
As before, we prove that, for $y\in\partial B_\e$
$$
\frac{\partial\psi}{\partial\nu}(y)= O\left(\frac{||\Pi v_\e^o||}{\l_\e^{1/2}d_\e^2}\right)\quad\mbox{and}\quad \frac{\partial\D\psi}{\partial\nu}(y)= O\left(\frac{||\Pi v_\e^o||}{\l_\e^{1/2}d_\e^4}\right).
$$
Therefore
\begin{eqnarray*}
\int_{B_\e} \d_\e^{8-\e}\Pi v_\e^o w_\e \leq \frac{c||\Pi v_\e^o||}{\l_\e^{1/2}d_\e^4}\left( \frac{1}{\d_\e^{3/2}}\int_{\partial B_\e}|w_\e| + \frac{1}{\d_\e^{7/2}}\int_{\partial B_\e}|\D w_\e|\right)
\leq \frac{c||v_\e|| ||\Pi v_\e^o||}{(\l_\e d_\e)^{1/2}}.
\end{eqnarray*}
The proof of the lemma is completed.
\end{pf}
\begin{lem}\label{l:28}
For $\e$ small, we have
\begin{align*}
i)&\quad \int_{B_\e}\d_\e^{7-\e}v_\e \frac{1}{\l_\e}\frac{\partial \d_\e}{\partial a_k}= O\left( \frac{||\Pi v_\e^o||}{\l_\e^{1/2}}+ \frac{||v_\e||}{\l_\e d_\e^{1/2}}\right),\\
ii)&\quad \int_{B_\e}\d_\e^{7-\e}\th_\e v_\e \frac{1}{\l_\e}\frac{\partial \d_\e}{\partial a_k}= O\left( \frac{||\Pi v_\e^o||}{\l_\e d_\e}+ \frac{||v_\e||}{(\l_\e d_\e)^{3/2}}\right).
\end{align*} 
\end{lem}
\begin{pf}
 Claim $i)$ can be proved in the same way as Lemma \ref{l:26}, so we omit its proof. Claim $ii)$ follows from Proposition \ref{p:21} and claim $i)$.
\end{pf}

Let us now compute the contribution of the following integral which involves $v_\e^2$.
\begin{lem}\label{l:29}
Form $\e$ small, we have
$$
 \int_{B_\e}\d_\e^{7-\e}v_\e^2 \frac{1}{\l_\e}\frac{\partial \d_\e}{\partial a_k}= O\left( ||\Pi v_\e^o|| ||v_\e|| + \frac{||v_\e||^2}{(\l_\e d_\e)^{1/2}}\right).
$$
\end{lem}
\begin{pf}
Using \eqref{10} and the fact that the even part of $v_\e^2$ has no contribution to the integrals, we obtain
$$
 \int_{B_\e}\d_\e^{7-\e}v_\e^2 \frac{1}{\l_\e}\frac{\partial \d_\e}{\partial a_k}=  \int_{B_\e}\d_\e^{7-\e}\frac{1}{\l_\e}\frac{\partial \d_\e}{\partial a_k}(2v_\e-w_\e)w_\e + O\left(||\Pi v_\e^o|| ||v_\e||\right).
$$
Let $\Psi$ be the solution of
$$
\D^2 \Psi =  \d_\e^{7-\e}\frac{1}{\l_\e}\frac{\partial \d_\e}{\partial a_k}(2v_\e-w_\e)\quad\mbox{in}\quad B_\e;\quad \D\Psi=\Psi=0\quad\mbox{on} \quad\partial B_\e.
$$
Thus, as in the proof of Lemma \ref{l:26}, we obtain for $y\in\partial B_\e$
$$
\frac{\partial\Psi}{\partial\nu}(y)=O\left(\frac{||v_\e||}{\l_\e^{1/2} d_\e^2}\right),\quad \frac{\partial\D\Psi}{\partial\nu}(y)=O\left(\frac{||v_\e||}{\l_\e^{1/2} d_\e^4}\right)
$$
and therefore
$$
\int_{B_\e}\d_\e^{7-\e}\frac{1}{\l_\e}\frac{\partial \d_\e}{\partial a_k}(2v_\e-w_\e)w_\e= O\left(\frac{||v_\e||^2}{(\l_\e d_\e)^{1/2}}\right).
$$
Thus our lemma follows.
\end{pf}

Next we are going to estimate the integrals involving the odd part of $v_\e$ with respect to $(x-a_\e)_k$, for $1\leq k\leq 5$.
\begin{lem}\label{l:210}
For $\e$ small, we have
$$
\int_{B_\e}u_\e^{9-\e}\Pi v_\e^o= 9\int_{B_\e}\d_\e^{8}(\Pi v_\e^o)^2 + o(||\Pi v_\e^o||^2) + O\left( ||\Pi v_\e^o||\left(\e^{3/2} + \frac{1}{(\l_\e d_\e)^{3/2}}\right)\right).
$$
\end{lem}
\begin{pf}
We have
\begin{align}\label{20}
\int_{B_\e}u_\e^{9-\e}\Pi v_\e^o &= \a_\e^{9-\e}\int_{B_\e}P\d_\e^{9-\e}\Pi v_\e^o + (9-\e) \a_\e^{8-\e}\int_{B_\e}P\d_\e^{8-\e}v_\e \Pi v_\e^o\notag\\
&+ O\left(\int_{B_\e}P\d_\e^{7-\e}|v_\e|^2|\Pi v_\e^o|+ \int_{B_\e}|v_\e|^{9-\e}|\Pi v_\e^o|\right)\notag\\
&= \a_\e^{9-\e}\int_{B_\e}P\d_\e^{9-\e}\Pi v_\e^o + (9-\e) \a_\e^{8-\e}\int_{B_\e}P\d_\e^{8-\e}v_\e \Pi v_\e^o + O(||v_\e||^2||\Pi v_\e^o||).
\end{align}
We estimate the two integrals on the right-hand side in \eqref{20}. First, using Proposition \ref{p:21} and the Holder inequality, we have
\begin{align*}
\int_{B_\e}P\d_\e^{8-\e}v_\e \Pi v_\e^o&= \int_{B_\e}\d_\e^{8-\e}v_\e \Pi v_\e^o + O\left(\frac{||v_\e||||\Pi v_\e^o|| }{\l_\e d_\e}\right)\\
&= \int_{B_\e} \d_\e^{8-\e} (\Pi v_\e^o)^2 + \int_{B_\e} \d_\e^{8-\e} \Pi v_\e^o w_\e,
\end{align*}
where we have used in the last equality the evenness of $\d_\e$ and $\Pi v_\e^e$ and the oddness of $\Pi v_\e^o$.\\
By Lemmas \ref{l:23} and \ref{l:27} we obtain
\be\label{21}
\int_{B_\e}P\d_\e^{8-\e}v_\e \Pi v_\e^o= \int_{B_\e} \d_\e^{8} (\Pi v_\e^o)^2  + O\left( \frac{||v_\e|| ||\Pi v_\e^o||}{(\l_\e d_\e)^{1/2}}\right).
\ee
Secondly, we write
$$
\int_{B_\e}P\d_\e^{9-\e}\Pi v_\e^o= \int_{B_\e}\d_\e^{9-\e}\Pi v_\e^o - (9-\e)\int_{B_\e}\d_\e^{8-\e}\th_\e\Pi v_\e^o + O\left(\int_{B_\e}\d_\e^{7-\e}\th_\e^2|\Pi v_\e^o|\right).
$$
Thus, using the evenness of $\d_\e$, the oddness of $\Pi v_\e^o$ and Holder inequality, we obtain
\be\label{22}
\int_{B_\e}P\d_\e^{9-\e} \Pi v_\e^o= O\left(\frac{||\Pi v_\e^o||}{(\l_\e d_\e)^2}\right).
\ee
Using \eqref{20}, \eqref{21}, \eqref{22} and Propositions \ref{p:24} and \ref{p:26}, we easily derive our lemma.   
\end{pf}
\begin{lem}\label{l:211}
For $\e$ small, we have
$$
||\Pi v_\e^o || = O\left(\e^{3/2} + \frac{1}{(\l_\e d_\e)^{3/2}}\right).
$$
\end{lem}
\begin{pf}
We write
 \be\label{23}
\Pi  v_\e^{o}=\tilde{\Pi v}_\e^o + \a \Pi\d_\e +
\b\l_\e\frac{\partial \Pi\d_\e}{\partial\l} + \sum_{r=1}^5\g_r
\frac{1}{\l_\e}\frac{\partial \Pi\d_\e}{\partial a_r} 
\ee
 with
$$
\langle \tilde{\Pi v}_\e^o, \Pi\d_\e\rangle = \langle \tilde{\Pi v}_\e^o,
\frac{\partial \Pi\d_\e}{\partial\l}\rangle = \langle \tilde{\Pi v}_\e^o,
\frac{\partial \Pi\d_\e}{\partial a_r}\rangle =0 \mbox{ for each } r\in
\{1,2,3,4,5\}.
$$
Taking the scalar product in $H^2\cap H^1_0(B_\e)$ of \eqref{23} with
$\Pi\d_\e$, $\l_\e\partial \Pi\d_\e / \partial\l$,
$\l_\e^{-1} \partial \Pi\d_\e / \partial a_r$, $1\leq r
\leq 5$, provides us with the following invertible linear system
in $\a$, $\b$, $\g_r$ (with  $1\leq r\leq 5$)
$$
(S)\,\,
\begin{cases}
\langle \Pi\d_\e, \Pi v_\e^o\rangle= \a(C'+o(1))+ \b\langle \Pi\d_\e,
\l_\e\frac{\partial \Pi\d_\e}{\partial\l}\rangle +\sum_{r=1}^5\g_r
\langle \Pi\d_\e, \frac{1}{\l_\e}\frac{\partial \Pi\d_\e}{\partial a_r}\rangle\\
\langle\l_\e\frac{\partial \Pi\d_\e}{\partial\l}, \Pi v_\e^{o}\rangle=
\a \langle \Pi\d_\e, \l_\e\frac{\partial \Pi\d_\e}{\partial\l}\rangle +
\b(C'' +o(1)) + \sum_{r=1}^5\g_r \langle \l_\e\frac{\partial
\Pi\d_\e}{\partial\l},
\frac{1}{\l_\e}\frac{\partial \Pi\d_\e}{\partial a_r}\rangle\\
\langle\frac{1}{\l_\e}\frac{\partial \Pi\d_\e}{\partial a_k},
\Pi v_\e^{o}\rangle= \a \langle \Pi\d_\e, \frac{1}{\l_\e}\frac{\partial
\Pi\d_\e}{\partial a_k}\rangle+ \b\langle \l_\e\frac{\partial
\Pi\d_\e}{\partial\l}, \frac{1}{\l_\e}\frac{\partial
\Pi\d_\e}{\partial a_k}\rangle + \sum_{r=1}^5\g_r \langle
\frac{1}{\l_\e}\frac{\partial \Pi\d_\e}{\partial a_k},
\frac{1}{\l_\e}\frac{\partial \Pi\d_\e}{\partial a_r}\rangle.
\end{cases}
$$
Observe that
$$
\langle \Pi\d_\e, \l_\e\frac{\partial \Pi\d_\e}{\partial\l}\rangle
=O\left(\frac{1}{\l_\e d_\e}\right);\quad \langle \l_\e\frac{\partial
\Pi\d_\e}{\partial\l}, \frac{1}{\l_\e}\frac{\partial \Pi\d_\e}
{\partial a_r}\rangle=O\left(\frac{1}{(\l_\e d_\e)^2}\right);
$$
$$
 \langle \Pi\d_\e, \frac{1}{\l_\e}\frac{\partial
\Pi\d_\e}{\partial a_r}\rangle=
O\left(\frac{1}{(\l_\e d_\e)^2}\right);\quad \langle
\frac{1}{\l_\e}\frac{\partial \Pi\d_\e}{\partial a_k},
\frac{1}{\l_\e} \frac{\partial
\Pi\d_\e}{\partial a_r}\rangle=(C'''+o(1))\d_{kr}+
O\left(\frac{1}{(\l_\e d_\e)^2}\right),
$$
where $\d_{kr}$ denotes the Kronecker symbol.\\
 Now, because of the evenness of $\d_\e$ and the oddness
of $\Pi v_\e^{o}$ with respect to $(x-a_\e)_k$ we obtain
 \be\label{24}
\langle \Pi\d_\e, \Pi v_\e^{o}\rangle = \int_{B_\e}\D \Pi\d_\e . \D \Pi v_\e^{o} = \int_{B_\e} \d_\e
^9 \Pi v_\e^{o} = 0.
 \ee
 In the same way we have
$$
\langle\l_\e\frac{\partial \Pi\d_\e}{\partial\l}, \Pi v_\e^{o}\rangle=
\langle\frac{1}{\l_\e}\frac{\partial \Pi\d_\e}{\partial a_r},
\Pi v_\e^{o}\rangle=0\quad\mbox{ for each } r\ne k.
$$
We also have
 \begin{align}\label{25}
 \langle \frac{1}{\l_\e} \frac{\partial\Pi\d_\e}{\partial a_k}, \Pi v_\e^{o}\rangle& =
\int_{B_\e}\D\left( \frac{1}{\l_\e} \frac{\partial\Pi\d_\e}{\partial a_k}\right) . \D\left(v_\e-\Pi v_\e^{e}-w_\e\right)\notag\\
& =\int_{B_\e}\D\left( \frac{1}{\l_\e} \frac{\partial\Pi\d_\e}{\partial a_k}\right) . \D v_\e -
\int_{B_\e}\D\left( \frac{1}{\l_\e} \frac{\partial\Pi\d_\e}{\partial a_k}\right) . \D w_\e,
\end{align}
where we have used in the last equality the fact that $\Pi v_\e^e$ is even with respect to $(x-a_\e)_k$.\\
Using \eqref{m:2} and Holder inequality, we obtain
\be\label{26}
\int_{B_\e}\D\left( \frac{1}{\l_\e} \frac{\partial\Pi\d_\e}{\partial a_k}\right) . \D v_\e \leq c||v_\e|| \left( \int_{\O\setminus B_\e} \bigg|\D \frac{1}{\l_\e}\frac{\partial \d_\e}{\partial a_k}\bigg|^2\right)^{1/2} = O\left(\frac{||v_\e||}{(\l_\e d_\e)^{3/2}}\right).
\ee
\eqref{26} and Lemma \ref{l:27} imply that
\be\label{27}
 \langle \frac{1}{\l_\e} \frac{\partial\Pi\d_\e}{\partial a_k}, \Pi v_\e^{o}\rangle
 = O\left(\frac{||v_\e||}{(\l_\e d_\e)^{3/2}}\right).
\ee
Inverting the linear system $(S)$, we deduce from the above estimates
\be\label{28}
\a=O\left(\frac{||v_\e||}{(\l_\e d_\e)^{\frac{7}{2}}}\right),\,
\b=O\left(\frac{||v_\e||}{(\l_\e d_\e)^{\frac{7}{2}}}\right), \,
\g_k=O\left(\frac{||v_\e||}{(\l_\e d_\e)^{\frac{3}{2}}}\right), \,
\g_r=O\left(\frac{||v_\e||}{(\l_\e d_\e)^{\frac{7}{2}}}\right), \,
r\ne k.
 \ee
 This implies through \eqref{23}
 \be\label{29}
||\Pi v_\e^{o}-\tilde{\Pi v}_\e^o||=
O\left(\frac{||v_\e||}{(\l_\e d_\e)^{3/2}}\right),\,\,\,
||\Pi v_\e^{o}||^2=||\tilde{\Pi v}_\e^o||^2+
O\left(\frac{||v_\e||^2}{(\l_\e d_\e)^3}\right).
 \ee
We now turn to the last step, which consists in estimating
$||\tilde{\Pi v}_\e^o||$. Since $u_\e$ is a solution of $(P_\e)$, we have
\be\label{30}
\int_{B_\e}\D^2 u_\e \Pi v_\e^o = \int_{B_\e}u_\e^{9-\e}\Pi v_\e^o.
\ee
Because of the evenness of $\d_\e$ and the oddness of $\Pi v_\e^o$ with respect to $(x-a_\e)_k$, \eqref{30} becomes
\be\label{31}
\mid\mid \Pi v_\e^o\mid\mid^2 = \int_{B_\e}u_\e^{9-\e}\Pi v_\e^o.
\ee
By \eqref{29}, \eqref{31} and Lemma \ref{l:210}, we obtain 
\be\label{k:2}
\mid\mid \tilde{\Pi v_\e^o}\mid\mid^2 -9 \int_{B_\e}\d_\e^8(\Pi v_\e^o)^2+o(||\tilde{\Pi v_\e^o}||^2)= O\left( \e^3+\frac{1}{(\l_\e d_\e)^3}\right) .
\ee
Using now \eqref{k:2} and the fact that the quadratic form 
$$
v\mapsto \int_{B_\e}|\D v|^2 - 9\int_{B_\e} \d_i^8v^2
$$
is positive definite (see \cite{BE}) on the subset
$\left[\mbox{Span}\,\left(\Pi\d_\e, \frac{\partial
\Pi\d_\e}{\partial\l}, \frac{\partial
\Pi\d_\e}{\partial a_k}\,\,1\leq k\leq
5\right)\right]^\bot_{H^2\cap H^1_0(B_\e)}$, we obtain
\be\label{32}
||\tilde{\Pi v_\e^o}||\leq C\left(\frac{1}{(\l_\e d_\e)^{3/2}}+\e^{3/2}\right).
\ee
Our lemma follows from \eqref{29} and \eqref{32}.
\end{pf}

Before ending this section, let us prove the following estimate which will be needed later.
\begin{lem}\label{l:212}
For $\e$ small, we have
$$
\langle \frac{\partial^2 P\d_\e}{\partial\l\partial a_k}, v_\e\rangle = O\left(\frac{1}{(\l_\e d_\e)^{3/2}} + \e^{3/2}\right).
$$
\end{lem}
\begin{pf}
We have
\begin{align}\label{33}
\int_{\O} \D\left(\frac{\partial^2 P\d_\e}{\partial\l\partial a_k}\right)\D v_\e&= \int_{B_\e} \D^2\left(\frac{\partial^2 P\d_\e}{\partial\l\partial a_k}\right) v_\e + O\left(\frac{|| v_\e||}{(\l_\e d_\e)^{9/2}}\right)\notag\\
&= \int_{B_\e} \D^2\left(\frac{\partial^2 P\d_\e}{\partial\l\partial a_k}\right) \Pi v_\e^o  +\int_{B_\e} \D^2\left(\frac{\partial^2 P\d_\e}{\partial\l\partial a_k}\right) w_\e +  O\left(\frac{||v_\e||}{(\l_\e d_\e)^{9/2}}\right).
\end{align}
For the first integral on the right-hand side in \eqref{33}, we have
\be\label{34}
 \int_{B_\e} \D^2\left(\frac{\partial^2 P\d_\e}{\partial\l\partial a_k}\right) \Pi v_\e^o = O\left(||\Pi v_\e^o||\right)= 
O\left(\frac{1}{(\l_\e d_\e)^{3/2}} + \e^{3/2}\right),
\ee
where we have used in the last equality Lemma \ref{l:211}.\\
Now let $\psi_4$ be the solution of
$$
\D^2 \psi_4 = \D^2\left(\frac{\partial^2 P\d_\e}{\partial\l\partial a_k}\right) \quad\mbox{in}\quad B_\e,\quad  \D\psi_4=\psi_4=0\quad\mbox{on}\quad \partial B_\e.
$$
Thus, as in the proof of Lemma \ref{l:26}, we obtain for $y\in \partial B_\e$
$$
\frac{\partial\psi_4}{\partial\nu}(y) = O\left(\frac{1}{\l_\e^{1/2}d_\e^2}\right), \quad \frac{\partial\D\psi_4}{\partial\nu}(y) = O\left(\frac{1}{\l_\e^{1/2}d_\e^4}\right)
$$
and therefore
\be\label{35}
\int_{B_\e} \D^2\left(\frac{\partial^2 P\d_\e}{\partial\l\partial a_k}\right) w_\e = O\left(\frac{||v_\e||}{(\l_\e d_\e)^{1/2}}\right).
\ee
From \eqref{33}, \eqref{34}, \eqref{35} and Proposition \ref{p:26}, we easily deduce our lemma.
\end{pf}
\section{Proof of Theorems}
\mbox{}
Let us start by proving the following crucial result:
\begin{pro}\label{p:31}
 For $u_\e=\a_\e P\d_{a_\e,\l_\e}+v_\e$ solution of $(P_{\e})$ with $\l_\e^\e =1 + o(1)$ as $\e$ goes to zero, we have the following estimates
$$
c_2\e+O(\e^2)-c_1\frac{H(a_\e,a_\e)}{\l_\e}+o\left(\frac{1
}{\l_\e d_\e}\right)=0, \leqno{(a)}
$$
$$
\frac{c_3}{\l_\e ^2}\frac{\partial H(a_\e,a_\e) }{\partial a} + o\left(\frac{1}{(\l_\e d_\e)^2}\right) + O\left(\e^{5/2} + \frac{\e\log (\l_\e)}{(\l_\e d_\e)^2} + \frac{1}{(\l_\e d_\e)^{5/2}}\right)=0, \leqno{(b)}
$$
where $c_1$, $c_2$ are the constants defined in Theorem \ref{t:11},
and where $c_3>0$.
\end{pro}
\begin{pf}
Since  claim $(a)$ was proved in \cite{BE2}, we only need to prove claim $(b)$.
Multiplying the equation $(P_{\e})$ by $\frac{1}{\l_\e}\frac {\partial
P\d_\e}{\partial a_k}$ and integrating on $\O$, we obtain for $1\leq k\leq 5$
\begin{align}\label{e:31}
0&=\int _{\O}\D^2 u_\e \frac{1}{\l _\e} \frac{\partial P\d _\e}{\partial a_k}
- \int _{\O}
u_\e ^{9-\e}\frac{1}{\l _\e} \frac{\partial P\d _\e}{\partial a_k}\nonumber\\
&= \a _\e \int _{\O} \d _\e ^9 \frac{1}{\l _\e} \frac{\partial P\d
_\e}{\partial a_k}- \int _{\O} \biggl[ (\a _\e P\d _\e )^{9-\e} +
(9-\e ) (\a _\e P\d _\e )^{8-\e}v_\e \nonumber\\
& + \frac{(9-\e)(8-\e)}{2}(\a_\e P\d_\e)^{7-\e} v_\e^2 \biggr]\frac{1}{\l _\e} \frac{\partial P\d _\e}{\partial a_k} + O\left( ||v_\e||^3\right).
\end{align}
We  estimate each term on the right-hand side in \eqref{e:31}.
First, by Proposition \ref{p:21} and the Holder inequality, we have 
\be \label{k:4}
\int_{\O} P\d_\e^{7-\e} v_\e^2 \frac{1}{\l _\e} \frac{\partial P\d _\e}{\partial a_k}=\int_{\O} \d_\e^{7-\e} v_\e^2 \frac{1}{\l _\e} \frac{\partial \d _\e}{\partial a_k} + O\left(\frac{||v_\e||^2}{\l_\e d_\e}\right).
\ee
Secondly, we compute
\begin{align}\label{e:32}
\int _{\O} P\d _\e ^{8-\e}& v_\e \frac{1}{\l _\e} \frac{\partial P\d
_\e}{\partial a_k}=\int _{\O} \d _\e ^{8-\e} v_\e \frac{1}{\l _\e} \frac{\partial P\d_\e}{\partial a_k} + (8-\e)\int _{\O} \d _\e ^{7-\e}\th_\e v_\e \frac{1}{\l _\e} \frac{\partial P\d_\e}{\partial a_k}+ O\left( \int _{\O} \d _\e ^{7-\e}\th_\e^2 |v_\e|\right)\notag\\
& = \int _{\O} \d _\e ^{8-\e} v_\e \frac{1}{\l _\e} \frac{\partial \d_\e}{\partial a_k} + O\left(\int _{\O} \d _\e ^{8-\e} |v_\e| \mid\frac{1}{\l _\e} \frac{\partial \th_\e}{\partial a_k}\mid\right) + (8-\e)\int _{\O} \d _\e ^{7-\e}\th_\e v_\e \frac{1}{\l _\e} \frac{\partial \d_\e}{\partial a_k}\notag\\
&+ O\left(\int _{\O} \d _\e ^{7-\e} \th_\e|v_\e| \mid\frac{1}{\l _\e} \frac{\partial \th_\e}{\partial a_k}\mid\right) +  O\left(\int _{\O} \d _\e ^{7-\e}\th_\e^2 |v_\e|\right).
\end{align}
By Proposition \ref{p:21} and the Holder inequality, we obtain
\be\label{e:33}
\int _{\O} \d _\e ^{7-\e}\th_\e^2 |v_\e|=O\left(\frac{||v_\e||}{(\l_\e d_\e)^2}\right),\quad  \int _{\O} \d _\e ^{7-\e} \th_\e|v_\e| \mid\frac{1}{\l _\e} \frac{\partial \th_\e}{\partial a_k}\mid =O\left(\frac{||v_\e||}{(\l_\e d_\e)^3,}\right),
\ee
\be\label{e:34}
\int _{\O} \d _\e ^{8-\e} |v_\e| \mid\frac{1}{\l _\e} \frac{\partial \th_\e}{\partial a_k}\mid= O\left(\frac{||v_\e||}{(\l_\e d_\e)^2}\right).
\ee
We also have by Proposition \ref{p:24}
\begin{align*}
 \int _{\O} \d _\e ^{8-\e} v_\e \frac{1}{\l _\e} \frac{\partial \d_\e}{\partial a_k}&=  \int _{\O} \d _\e ^{8}\left(\d_\e^{-\e}-\frac{c_0^{-\e}}{\l_\e^{\e/2}}\right) v_\e \frac{1}{\l _\e} \frac{\partial \d_\e}{\partial a_k}\\
&= \int _{B_\e} \d _\e ^{8}\left(\d_\e^{-\e}-\frac{c_0^{-\e}}{\l_\e^{\e/2}}\right) v_\e \frac{1}{\l _\e} \frac{\partial \d_\e}{\partial a_k}+ O\left(\frac{||v_\e||}{(\l_\e d_\e)^{9/2}}\right).
\end{align*}
Using \eqref{10}, Lemma \ref{l:23} and the Holder inequality, we derive that
\begin{align}\label{e:35}
\int _{\O} \d _\e ^{8-\e} v_\e \frac{1}{\l _\e} \frac{\partial \d_\e}{\partial a_k}&=\int _{B_\e} \d _\e ^{8}\left(\d_\e^{-\e}-\frac{c_0^{-\e}}{\l_\e^{\e/2}}\right) \frac{1}{\l _\e} \frac{\partial \d_\e}{\partial a_k}w_\e +O\left(\e ||\Pi v_\e^o|| + \frac{||v_\e||}{(\l_\e d_\e)^{9/2}}\right)\notag\\
&=O\left(\frac{\e ||v_\e||}{(\l_\e d_\e)^{1/2}} + \e ||\Pi v_\e^o|| + \frac{||v_\e||}{(\l_\e d_\e)^{9/2}}\right),
\end{align}
where we have used Lemma \ref{l:26} in the last equality.\\
Using \eqref{k:4}--\eqref{e:35}, Lemmas \ref{l:25}, \ref{l:28}, \ref{l:29}, Proposition \ref{p:24} and the fact that $\l_\e ^\e= 1+O(\e \log\l_\e)$, we easily derive our result.
\end{pf}

We are now able to prove Theorem \ref{t:11}.\\
\begin{pfn}{Theorem \ref{t:11}}
Let $(u_\e)$ be a solution of $(P_{\e})$ which satisfies $(H)$.
Then, using Proposition \ref{p:24}, $u_\e=\a_\e
P\d_{a_\e,\l_\e}+v_\e$ with $\a_\e\to 1$, $\l_\e ^\e\to 1$,  $\l_\e d(a_\e, \partial
\O) \to \infty$, $v_\e$ satisfies $(V_0)$ and $|| v_\e||\to 0$.
Now, using claim $(a)$ of Proposition \ref{p:31}, we derive that
\begin{eqnarray}\label{e:39}
\e=\frac{c_1}{c_2}\frac{H(a_\e,a_\e)}{\l_\e
}+o\left(\frac{1}{(\l_\e d_\e)^{n-4}}\right)=O\left(\frac{1}{\l_\e
d_\e}\right).
\end{eqnarray}
 Therefore, it follows from claim (b) of Proposition \ref{p:31} and Lemma \ref{l:211} that 
\begin{eqnarray}\label{e:310}
 \frac{\partial H(a_\e,a_\e) }{\partial a}=o\left(\frac{1}{d_\e
 ^2}\right).
 \end{eqnarray}
Using \eqref{e:310} and the fact that for $a$ near the boundary $\frac{\partial H}{\partial a} (a_\e,a_\e)\sim cd(a_\e,\partial\O)^{-2}$, we derive that $a_\e$ is away from the boundary and it converges to a
critical point $x_0$ of $\var$.\\
 Finally, using \eqref{e:39}, we obtain
 $$
\e \l_\e \to \frac{c_1}{c_2}\var (x_0) \mbox{ as } \e \to 0.
$$
By Proposition \ref{p:24}, we have
 \begin{eqnarray}\label{e:311}
||u_\e||_{L^\infty}^2\sim c_0^2\l_\e \quad \mbox{as}
\quad \e\to 0.
\end{eqnarray}
 This concludes the proof of Theorem \ref{t:11}.
\end{pfn}

The sequel of this section is devoted to the proof of Theorem \ref{t:12}.

\begin{pfn}{ Theorem \ref{t:12}}
Let $x_0$ be a nondegenerate critical point of $\var$. It is easy
to see that $d(a, \partial \O) > d_0 > 0$ for $a$ near $x_0$. We will take a
function $u=\a P\d_{(a,\l)} +v$ where $(\a-\a_0)$ is very small,
$\l$ is large enough, $||v||$ is very small, $a$ is close to $x_0$
and $\a_0=S^{-5/8}$ and we will prove that we can choose the
variables $(\a, \l, a, v)$ so that $u$ is a critical point of
$J_\e$ with $||u||=1$. Here $J_\e$ denotes the functional corresponding to problem $(P_{\e})$ defined by
$$
J_\e (u)=\left(\int_\O|\D u|^2\right) \left(\int_\O |u|^{10-\e}\right)^{-2/(10-\e)}.
$$
Let
\begin{align*}
M_\e =\{ (\a, \l, a, v)\in \R^*_+\times \R^*_+ \times \O  \times
\mathcal{H}(\O)/ & \,  |\a-\a_0|< \nu_0,\  d_a > d_0,\  \l >
\nu_0^{-1}, \\
 & \e \log\l < \nu_0,\  ||v||<\nu_0 \mbox{ and } v\in
E_{(a,\l)}\},
\end{align*}
where $\nu_0$ and $d_0$ are two suitable positive constants and where $d_a=d(a,\partial\O )$.\\
 Let us define
the functional
\begin{eqnarray*}
K_\e :  M_\e \to \R, \quad  K_\e (\a, a, \l, v) = J_\e(\a P\d_{(a,
\l)} +v).
\end{eqnarray*}
It is known  that $(\a, \l, a, v)$ is a critical point of $K_\e$
if and only if $u=\a P\d_{(a,\l)} + v$ is a critical point
of $J_\e$ on $E$. So this fact allows us to look for
critical points of $J_\e$ by successive optimizations with
respect to the different parameters on $M_\e$.\\
First, we know that (see \cite{BE2}) the following problem
$$\min \{ J_\e(\a P\d_{(a, \l)} +v),
\ v \,\,\mbox{satisfying } (V_0)\,\, \mbox{ and } ||v||<\nu_0\}
$$
is achieved by a unique function $\ov{v}$ which satisfies the
estimate of Proposition \ref{p:26}. This implies that there exist
$A$, $B$ and $C_i$'s such that
\begin{align}\label{ou1}
\frac{\partial K_\e}{\partial v}(\a, \l, a, \ov{v})&=\n J_\e (\a
P\d_{(a, \l)} +\ov{v})\notag\\
&=A  P\d_{(a, \l)}+B \frac{\partial
}{\partial \l}P\d_{(a, \l)}+\sum_{i=1}^5C_i \frac{\partial
}{\partial a_i}P\d_{(a, \l)},
\end{align}
where $a_i$ is the $i^{th}$ component of $a$.\\
According to \cite{BE2}, we have that
$$
A =O\left(\e \log \l +|\b|+\frac{1}{\l }\right), \quad B
=O\left(\l \e + 1\right),\quad C_j
=O\left(\frac{\e ^2}{ \l} +\frac{1}{\l ^{3}}\right).
$$
 To find critical points of $K_\e$, we have to solve the following system
$$
(E_1)\qquad
\begin{cases}
\frac{\partial K_\e}{\partial \a} & =0\\
\frac{\partial K_\e}{\partial \l} & = B\langle \frac{\partial ^2
P\d}{\partial \l ^2}, \bar{v}\rangle+\sum_{i=1}^5 C_i
\langle\frac{\partial ^2 P\d}{\partial \l \partial a_i},
\bar{v}\rangle\\
\frac{\partial K_\e}{\partial a_j} & =B\langle\frac{\partial ^2
P\d}{\partial \l \partial a_j}, \bar{v}\rangle+\sum_{i=1}^5 C_i
\langle\frac{\partial ^2 P\d}{\partial a_i\partial a_j}, \bar{v}\rangle,
\mbox{ for each } j=1,...,5.
\end{cases}
$$
Observe that for $\psi=P\d_{(a, \l)}$,  $\partial
P\d_{(a, \l)}/\partial \l$, $ \partial P\d_{(a,
\l)}/\partial a_i$ with $i=1,...,5$ and
 for $u=\a P\d_{(a, \l)}+\ov{v}$, we have
$$
\langle\n J_\e(u), \psi\rangle = 2J_\e (u) \left( \a \langle P\d_{(a, \l)},
\psi\rangle-J_\e (u)^{5-\e/2}\int _\O |u|^{8-\e}u\psi\right).
$$
We also have (see \cite{BE2})
\begin{eqnarray}\label{e:*1}
J_\e(\a P\d_{(a,\l)}+\ov{v})=S+O\left( \e \log\l +\frac{1}{\l}\right),
\end{eqnarray}
\begin{eqnarray*}
\frac{\partial K_\e}{\partial \a}  =\langle\n J_\e(\a P\d +\ov{v}), P\d\rangle
  = 2J_\e(u) \left( \a S^{5/4} \left(1-\a ^{8} S^{5}\right)
 + O\left(\e \log \l +\frac{1}{\l}\right)\right),
 \end{eqnarray*}
and
\begin{align*}
\l \frac{\partial K_\e}{\partial \l} & =\langle\n J_\e(\a P\d
+\ov{v}),\l \frac{\partial P\d}{\partial \l}\rangle\\
 & = J_\e(u) \left( \a c_1\frac{H(a,a)}{\l}\left(1-2
\a ^{8} S^{5}\right) +c_2S^{5}\a ^9 \e +O\left(\e ^2 \log \l +\frac{\e \log \l}{\l}+\frac{1}{\l
^3} \right)\right).
 \end{align*}
Following  the proof of claim (b) of Proposition \ref{p:31},
we obtain, for each \\ $j=1,...,5$,
\begin{align*}
\frac{1}{\l}\frac{\partial K_\e}{\partial a_j} & =\langle\n J_\e(\a P\d +\ov{v}),\frac{1}{\l}\frac{\partial P\d}{\partial a_j}\rangle\\
&  = - \frac{c\a }{2\l^2}\frac{\partial H(a,a)}{\partial a}\left(1-2
\a ^8 S^{5}\right) + O\left( \e ^{5/2}  +\frac{\e \log
\l}{\l ^2}+\frac{1}{\l ^{5/2}}\right).
 \end{align*}
On the other hand, one can easily verify that
\be\label{e:R}
(i)\, \bigg|\bigg|\frac{\partial^2 P\d}{\partial \l^2}\bigg|\bigg|=O\left(\frac{1}{\l^2}\right), \,  \, (ii)\, \bigg|\bigg|\frac{\partial^2 P\d}{\partial a_i \partial a_j}\bigg|\bigg|=O(\l^2).
\ee
Now, we take the following change of variables:
$$
\a= \a_0 +\b, \quad a=x_0+\xi, \quad \frac{1}{\l ^{\frac{1}{2}}}
= \sqrt{\frac{c_2}{c_1}}\left( \frac{1}{\sqrt{H(x_0,x_0)}}+\rho\right)
\sqrt{\e}.
$$
Then, using estimates \eqref{e:R}, Lemma \ref{l:212}, Proposition \ref{p:26} and the fact that $x_0$ is a nondegenerate critical point
of $\var$, the system $(E_1)$ becomes
$$
(E_2) \qquad
\begin{cases}
\b & = O\left( \e |\log \e|+|\b|^2\right)\\
\rho & =O\left( \e |\log\e| +|\b|^2+|\xi|^2+\rho ^2\right)\\
\xi & =O\left( |\b|^2+|\xi|^2+\e ^{1/2}\right).
\end{cases}
$$
Thus Brower's fixed point theorem shows that the system $(E_2)$ has a solution $(\b_\e,
\rho_\e, \xi_\e)$ for $\e$ small enough such that
$$
\b_\e = O(\e |\log \e|), \quad \rho_\e =O(\e|\log\e|), \quad \xi_\e=O( \e^{1/2}).
$$
 By construction, the corresponding $u_\e$ is a critical point of
 $J_\e$ that is  $w_\e = J_\e(u_\e)^{(5-\e/2)/(8-\e)}u_\e$ satisfies
 \be\label{e:ew}
 \D ^2 w_\e = |w_\e|^{8 - \e}w_\e \mbox{ in }\O, \quad w_\e=\D
 w_\e=0 \mbox{ on } \partial \O
 \ee
with $|w_\e^-|_{L^{10}(\O)}$ very small, where $w_\e
^-=\max(0,-w_\e)$.\\
As in Proposition 4.1 of \cite{BEH2}, we prove that $w_\e^-=0$. Thus, since $w_\e$ is
a non-negative function which satisfies \eqref{e:ew}, the strong maximum
principle ensures that $w_\e > 0$ on $\O$ and then $w_\e$ is a
solution of $(P_{\e})$, which blows up at $x_0$ as $\e$ goes to
zero. This ends the proof of Theorem \ref{t:12}.
\end{pfn}

\end{document}